\documentclass[11pt,reqno]{amsart}
\usepackage[usenames]{color}
\usepackage{courier}
\usepackage{mathrsfs}
\usepackage{hyperref}
\usepackage{amssymb,amsmath,amsfonts,latexsym,amsthm}
\usepackage{graphicx}
\usepackage{psfrag}
\usepackage{float}
\usepackage{url}		

\usepackage{tikz}
\usetikzlibrary{positioning,arrows.meta}
\usetikzlibrary{calc}

\calclayout

\newcommand{\seqnum}[1]{\href{http://oeis.org/#1}{\underline{#1}}}

\newcount\tilpos

\newcommand{\Tiling}[2][4]{% \Tiling[<N>]{A,B,C,...}
\begin{tikzpicture}[x=0.7cm,y=0.7cm,baseline={(0,0.5)}, line join=round]
  \def\N{#1}

  \tikzset{
    board/.style={draw, line width=1.0pt},
    Astyle/.style={draw, line width=0.9pt}, % light green
    Bstyle/.style={draw, line width=0.9pt},  % light blue
    Cstyle/.style={draw, line width=0.9pt},   % light red
    lab/.style={font=\small}
  }

  \def\LenA{1}\def\LenB{1}\def\LenC{2}

  \def\TileA##1{%
    \draw[Astyle] (##1,0) rectangle ++(1,1);
    \node[lab] at ({##1+0.5},0.5) {$at$};
  }
  \def\TileB##1{%
    \draw[Bstyle] (##1,0) rectangle ++(1,1);
    \node[lab] at ({##1+0.5},0.5) {$b$};
  }
  \def\TileC##1{%
    \draw[Cstyle] (##1,0) rectangle ++(2,1);
    \node[lab] at ({##1+1.0},0.5) {$c$};
  }

  \tilpos=0
  \foreach \t in {#2} {%
    \edef\x{\number\tilpos}%
    \csname Tile\t\endcsname{\x}%
    \global\advance\tilpos by \csname Len\t\endcsname\relax
  }

  \draw[board] (0,0) rectangle (\N,1);
\end{tikzpicture}%
}

\newcommand{\SetAColor}[2]{%
  \expandafter\def\csname Acolor#1\endcsname{#2}%
}

\newcount\tilpos

\newcommand{\TilinC}[2][4]{% \TilinC[<N>]{A[1],B,C,...}
\begin{tikzpicture}[x=0.7cm,y=0.7cm,baseline={(0,0.5)}, line join=round]
  \def\N{#1}

  \tikzset{
    board/.style={draw, line width=1.0pt},
    Bstyle/.style={draw, line width=0.9pt, fill=white},
    Cstyle/.style={draw, line width=0.9pt, fill=white},
    lab/.style={font=\small}
  }

  \def\LenA{1}\def\LenB{1}\def\LenC{2}

  \def\sp{ }%
  \def\Trimmed{}%
  \def\TrimLead##1{\expandafter\TrimLeadAux##1\relax}%
  \def\TrimLeadAux##1##2\relax{%
    \ifx##1\sp
      \expandafter\TrimLeadAux##2\relax
    \else
      \def\Trimmed{##1##2}%
    \fi
  }%

  \def\SplitFirst##1##2\relax{\def\TileType{##1}\def\TileRest{##2}}%
  \def\LBrack{[}%
  \def\GetFirst##1##2\relax{\def\RestFirst{##1}\def\RestTail{##2}}%
  \def\StripBrackets[##1]\relax{\def\Aindex{##1}}%

  \def\TA{A}\def\TB{B}\def\TC{C}

  \tilpos=0
  \foreach \rawtile in {#2} {%
    \TrimLead{\rawtile}%
    \edef\x{\number\tilpos}%
    \expandafter\SplitFirst\Trimmed\relax

    \ifx\TileType\TA
      \def\Aindex{1}%
      \ifx\TileRest\empty\else
        \expandafter\GetFirst\TileRest\relax
        \ifx\RestFirst\LBrack
          \expandafter\StripBrackets\TileRest\relax
        \else
          \def\Aindex{\TileRest}%
        \fi
      \fi
      \ifcsname Acolor\Aindex\endcsname
        \edef\Afill{\csname Acolor\Aindex\endcsname}%
      \else
        \def\Afill{gray!20}%
      \fi
      \draw[draw, line width=0.9pt, fill=\Afill] (\x,0) rectangle ++(1,1);
      \node[lab] at ({\x+0.5},0.5) {$ax$};

    \else\ifx\TileType\TB
      \draw[Bstyle] (\x,0) rectangle ++(1,1);
      \node[lab] at ({\x+0.5},0.55) {$bx$};

    \else\ifx\TileType\TC
      \draw[Cstyle] (\x,0) rectangle ++(2,1);
      \node[lab] at ({\x+1.0},0.58) {$cx^2$};
    \fi\fi\fi

    \global\advance\tilpos by \csname Len\TileType\endcsname\relax
  }

  \draw[board] (0,0) rectangle (\N,1);
\end{tikzpicture}%
}

\SetAColor{1}{green!20}
\SetAColor{2}{yellow!20}
\SetAColor{3}{orange!25}

\makeatother

\newcommand{\N}{{\mathbb N}}

\theoremstyle{plain}
\numberwithin{equation}{section}
\newtheorem{thm}{Theorem}[section]
\newtheorem{corollary}{Corollary}[section]
\newtheorem{theorem}[thm]{Theorem}
\newtheorem{lemma}[thm]{Lemma}

\newtheorem{definition}[thm]{Definition}
\newtheorem{proposition}[thm]{Proposition}

\newtheorem{remark}{Remark}
\newtheorem{OpenP}{Open problem}

\title[On $p$-Adic Valuation Trees Arising from Laplace Transform of Fibonacci Polynomials]{On $p$-Adic Valuation Trees Arising from Laplace Transform of Fibonacci Polynomials}

\author[R. Fl\'orez]{Rigoberto Fl\'orez}
\address{Department of Mathematics and Data Analytics\\
		The Citadel\\
		Charleston, SC \\
		U.S.A.}
  \email{rflorez1@citadel.edu}
\author[R.  Higuita-Diaz]{Robinson Higuita-Diaz}
 \address{Escuela de Matem\'aticas y Estad\'istica, Universidad Pedag\'ogica y Tecnol\'ogica de Colombia, Tunja, Colombia.}
\email{ robinson.higuita@uptc.edu.co}
\author[J. L. Ram\'{\i}rez]{Jos\'e L. Ram\'{\i}rez}
\address{Departamento de Matem\'aticas,  Universidad Nacional de Colombia,  Bogot\'a, Colombia.}
\email{jlramirezr@unal.edu.co}

\date{\today}

\begin{document}

\maketitle

\begin{abstract} 
We study the $p$-adic valuation trees associated with the sequence \[ e_n=\int_0^\infty e^{-t}F_n(t)\,dt, \] obtained by evaluating the Laplace transforms of the Fibonacci polynomials at $1$. We answer an open question of Givens and Moll in the negative by constructing an explicit counterexample showing that even internal branches need not be simple. We then reformulate the Givens--Moll lifting criteria for the even and odd subsequences in terms of the roots and derivatives of two associated polynomials.  In particular, simple roots modulo $p$ give rise to simple internal branches.  Finally, we develop a recursive lifting procedure that determines the unique surviving residue class  at each level of a branch associated with a simple root and computes the $p$-adic valuations of the corresponding terms. 
\end{abstract}

\begin{center}
\emph{Dedicated to the memory of our friend and colleague Florian Luca.\footnote{Florian Luca had agreed to collaborate with us on this series of papers, but, unexpectedly, he passed away before he had the opportunity to communicate any results of his work to us. Although we were unable to pursue this endeavor with him, we remember with profound gratitude his friendship, generosity, encouragement, and the inspiration he provided. This paper is dedicated to his memory.}}
\end{center}

\section{Introduction}
Givens and Moll \cite{Givens} studied the sequence obtained by
evaluating the Laplace transforms of the Fibonacci polynomials at $1$.  More recently, Fl\'orez, Higuita, and Ram\'{\i}rez \cite{FlorezHiguitaRamirezHolonomic} extended this construction to other families of polynomials satisfying second-order recurrences, including the Lucas, Fermat, Morgan--Voyce, and Vieta polynomials.

Let $F_n(x)$ denote the $n$th Fibonacci polynomial, defined by 
\[ 
F_n(x)=xF_{n-1}(x)+F_{n-2}(x),\qquad n\geq2, 
\] 
with initial conditions $F_0(x)=0$ and $F_1(x)=1$. Givens and Moll considered the integer sequence $(e_n)_{n\geq1}$ defined by 
\[ 
e_n=\int_0^\infty e^{-t}F_n(t)\,dt. 
\] 
We refer to this sequence as the \emph{FiboLaplace sequence}. Its first terms are 
\[ 
1,\ 1,\ 3,\ 8,\ 31,\ 147,\ 853,\ 5824,\ 45741,\ 405845,\ 4012711,\ 43733976,\ldots. 
\] 
It appears in the OEIS as \cite[\seqnum{A003470}]{OEIS}

The FiboLaplace sequence admits the explicit representation \begin{equation}\label{eq:en-explicit} 
e_n= \sum_{i=0}^{\lfloor (n-1)/2\rfloor} \frac{(n-i-1)!}{i!}, \qquad n\geq1,
\end{equation} 
and satisfies the nonhomogeneous recurrence \begin{equation}\label{enRecurrenceRelation} 
e_m=(m-1)e_{m-1}-e_{m-2}+1-(-1)^m, \qquad m\geq3, 
\end{equation} 
with initial conditions $e_1=e_2=1$. Further properties and generalizations of this sequence can be found in \cite{FlorezHiguitaRamirezHolonomic,Givens}.

A central topic in \cite{Givens} is the study of the $p$-adic valuation of the terms of the FiboLaplace sequence. For a prime $p$, recall that $\nu_p(m)$ denotes the largest nonnegative integer $r$ such that $p^r\mid m$. For example, Givens and Moll obtained the following complete description of the $3$-adic valuation of $e_n$: 
\begin{equation}\label{theorem36moll} 
\nu_3(e_n)= 
\begin{cases} \nu_3(n), & \text{if } n\equiv0\pmod6,\\[1mm] 
    0, & \text{if } n\equiv1,2,4,5\pmod6,\\[1mm] 
    1, & \text{if } n\equiv3\pmod6. 
\end{cases} 
\end{equation}

Valuation trees have been used to study the $p$-adic valuations of
integer sequences; see, for instance, \cite{MollBook}. Examples
include Stirling numbers of the second kind \cite{BerrizbeitiaEtAl},
truncated exponential sums \cite{AmdeberhanCallanMoll}, Eulerian
numbers \cite{CastroGonzalezMedina}, and a sequence arising from a
rational integral \cite{SunMoll}.  Givens and Moll \cite{Givens}
constructed a tree structure for the FiboLaplace sequence that
organizes the congruence classes of indices according to the
$p$-adic valuations of the corresponding terms. A vertex of such a tree represents a congruence class of indices. It is called terminal when the valuation is constant on that class and nonterminal when further refinement is required. An internal branch is called simple when every level contains exactly one nonterminal vertex. For even indices, Givens and Moll established a lifting criterion that guarantees the simplicity of an internal branch. They subsequently asked whether every even internal branch must be simple.

The first objective of this paper is to answer this question. We show that the answer is negative by giving an explicit counterexample for $p=971$. We also observe that the analogous statement for odd internal branches already fails for $p=3$.

Our second objective is to provide an algebraic interpretation of the lifting criteria of Givens and Moll. We associate polynomials $\mathcal P_p(x)$ and $\mathcal Q_p(x)$ with the even and odd subsequences, respectively. Their roots modulo $p$ correspond to the internal roots of the associated valuation trees. Thus, the nonvanishing conditions in the Givens--Moll criteria are equivalent to the corresponding roots being simple modulo $p$. It follows that a simple root of $\mathcal P_p(x)$ or $\mathcal Q_p(x)$ gives rise to a simple even or odd internal branch, respectively.

Finally, we use the lifting congruences constructively. Along a simple branch, they determine the unique surviving residue class at each successive level and therefore provide a recursive method for computing the corresponding $p$-adic valuations. We illustrate this procedure in detail for $p=13$, where the four even internal branches are simple and their successive surviving residues can be computed one base-$13$ digit at a time. 

We first introduce the $p$-adic valuation trees associated with the FiboLaplace sequence and illustrate their structure through the case $p=13$. After recalling the lifting criteria of Givens and Moll, we present a counterexample for $p=971$ and then give an algebraic interpretation of the even and odd criteria in terms of the polynomials $\mathcal P_p(x)$ and $\mathcal Q_p(x)$. This interpretation leads to simple-root criteria and a partial converse. We conclude by developing a recursive procedure for following simple branches and computing the corresponding $p$-adic valuations, with a detailed illustration for $p=13$.

\section{Discussion on the Question of Givens and Moll}

In this section, we introduce the background needed to study the question of Givens and Moll \cite{Givens} concerning the $p$-adic valuations of the FiboLaplace sequence $(e_n)_{n\geq1}$. For a prime $p$ and a positive integer $m$, the \emph{$p$-adic valuation} of $m$, denoted by $\nu_p(m)$, is the largest nonnegative integer $r$ such that $p^r\mid m$. Thus, $\nu_p(m)$ is the exponent of $p$ in the prime factorization of $m$. For example, 
\[ 
\nu_2(100)=\nu_2(2^2\cdot 5^2)=2 \qquad\text{and}\qquad \nu_7(100)=0. 
\]

\begin{figure}[ht!]
    \centering
    \includegraphics[scale=1]{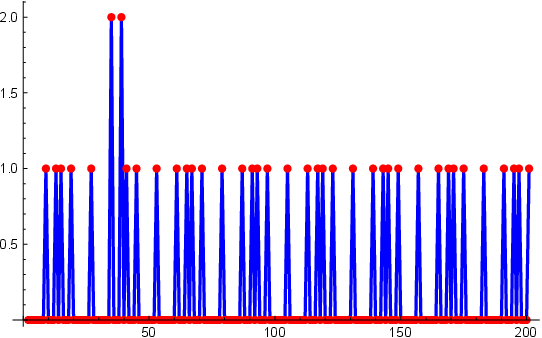}
    \caption{The $13$-adic valuation $\nu_{13}(e_n)$ for the first $200$ terms of the FiboLaplace sequence. }
    \label{fig1}
\end{figure}

To motivate the construction of the associated valuation trees, we begin with the case $p=13$. Figure~\ref{fig1} displays the values of $\nu_{13}(e_n)$ for $1\leq n\leq200$. Even in this relatively short range, a pattern is apparent: positive valuations occur at regularly spaced indices, and at some of these indices the valuation is greater than $1$.

\begin{figure}[ht!]
    \centering
    \includegraphics[scale=1]{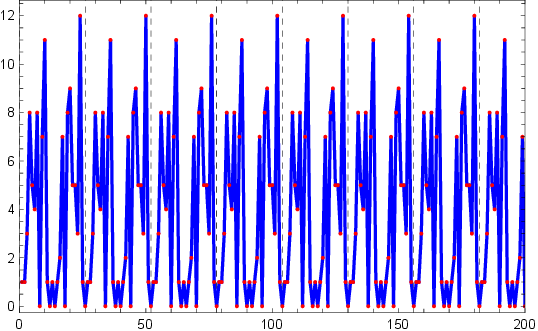}
    \caption{The residues $e_n\bmod 13$ for $1\leq n\leq 200$.
    The dashed vertical lines separate consecutive periods of length $26$.}
    \label{fig2}
\end{figure}

The first of these regularities is explained by a result of Givens and
Moll \cite{Givens}. They proved that, for every positive integer $M$, the sequence
$(e_n\bmod M)_{n\geq 1}$ is periodic with period $2M$. In particular,
\[
e_{n+26}\equiv e_n\pmod{13}.
\]
Consequently,
\[
13\mid e_n
\qquad\Longleftrightarrow\qquad
13\mid e_{n+26}.
\]
Thus divisibility by $13$ depends only on the residue class of $n$ modulo
$26$. Figure~\ref{fig2} displays the residues $e_n\bmod 13$ for
$1\leq n\leq 200$. The dashed vertical lines are spaced $26$ units apart,
making the periodicity particularly easy to see.

This suggests that we first divide the positive integers into the classes
\[
I_{13}(r)
=
\{n\geq 1:n\equiv r\pmod{26}\},
\qquad
1\leq r\leq 26.
\]
A computation over one period gives
\[
13\mid e_r
\qquad\Longleftrightarrow\qquad
r\in\{8,12,14,18,26\}.
\]
It follows immediately from periodicity that $\nu_{13}(e_n)=0$ 
for every $n$ outside the five classes
\begin{align*}
I_{13}(8)&=\{8,34,60,86,112,138,164,190,\dots\},\\
I_{13}(12)&=\{12,38,64,90,116,142,168,194,\dots\}, \\
I_{13}(14)&=\{14,40,66,92,118,144,170,196,\dots\}, \\
I_{13}(18)&=\{18,44,70,96,122,148,174,200,\dots\}, \\
I_{13}(26)&=\{26,52,78,104,130,156,182,208,\dots\}. 
\end{align*}

One might initially expect the residue class of $n$ modulo $26$ to determine not only whether $13\mid e_n$, but also the exact value of $\nu_{13}(e_n)$. However, this is not the case. For example, 
\[ 
8\equiv34\pmod{26}, 
\] 
whereas 
\[ 
\nu_{13}(e_8)=1 \qquad\text{and}\qquad \nu_{13}(e_{34})=2. 
\]
Thus, the residue class of $n$ modulo $26$ determines whether $\nu_{13}(e_n)>0$, but not its exact value. In particular, it does not distinguish between indices for which the valuation is $1$, $2$, or higher. This leads naturally to the following question:

\begin{center} 
\emph{What additional information is needed to distinguish the valuations within a class $I_{13}(r)$?} 
\end{center}

To detect whether $\nu_{13}(e_n)\geq2$, we must refine the congruence classes modulo $26$. By the periodicity theorem applied with $M=13^2$, the sequence \[ \bigl(e_n\bmod 13^2\bigr)_{n\geq1} \] has period \[ 2\cdot13^2=338. \] It is therefore natural to refine each residue class modulo $26$ into residue classes modulo $338$.

Hence a residue class modulo $26$ breaks into $13$ residue classes modulo
$338$. For example,
\[
I_{13}(8)  = \bigcup_{j=0}^{12} \left\{ n\geq 1: n\equiv 8+26j\pmod{338} \right\}.
\]
Thus the class $n\equiv 8\pmod{26}$ is refined into
\[
n\equiv 8,\ 34,\ 60,\ 86,\ 112,\ 138,\ 164,\ 190,\ 216,\ 242,\ 268,\ 294,\ 320 \pmod{338}.
\]

A computation modulo $13^2$ now reveals something more interesting.
Twelve of these classes contain only indices $n$ for which $\nu_{13}(e_n)=1$,
whereas the class $n\equiv 34\pmod{338}$
satisfies $\nu_{13}(e_n)\geq 2$. Thus, the original class $I_{13}(8)$ has exactly one subclass in which the
valuation can increase.

This refinement is naturally represented by a tree, which we denote by
$T_{13}(8)$. The root corresponds to the residue class $n\equiv 8\pmod{26}$. Its $13$ children, denoted by $v_{1,0},v_{1,1},\ldots,v_{1,12}$, correspond respectively to the classes
\[
n\equiv 8+26j\pmod{338},
\qquad 0\leq j\leq12.
\]
We place the value $1$ at a vertex when every index in the corresponding
class has $13$-adic valuation exactly $1$. If the valuation may be larger,
we write instead a lower bound. In the present example, $v_{1,1}$ is the unique vertex labeled $\geq2$, since it represents the class
\[
n\equiv34\pmod{338}.
\]
The other twelve vertices are terminal at this stage.

We now repeat the same procedure on the unique class that survives.
Since the sequence modulo $13^3$ has period $2\cdot13^3=4394$, the class $n\equiv34\pmod{338}$ splits into the $13$ classes
\[
n\equiv34+338j\pmod{4394},
\qquad 0\leq j\leq12.
\]
These form the second generation $
v_{2,0},v_{2,1},\ldots,v_{2,12}$ of the tree. Again, twelve classes have constant valuation $2$, while
exactly one continues $34+338=372$. 
Thus $v_{2,1}$ represents
\[
n\equiv372\pmod{4394}
\]
and is labeled $\geq3$.

\begin{figure}[ht]
\centering
\resizebox{0.9\textwidth}{!}{%
\begin{tikzpicture}[
    >={Stealth[length=1.4mm,width=0.9mm]},
    font=\scriptsize,
    root/.style={
        circle,
        draw,
        thick,
        minimum size=7mm,
        inner sep=1pt
    },
    term/.style={
        rectangle,
        draw,
        rounded corners=1pt,
        minimum width=6mm,
        minimum height=4.5mm,
        inner sep=0.8pt,
        align=center
    },
    cont/.style={
        rectangle,
        draw,
        rounded corners=1pt,
        thick,
        fill=red!15,
        minimum width=6mm,
        minimum height=4.5mm,
        inner sep=0.8pt,
        align=center
    },
    deadedge/.style={
        draw=gray,
        ->,
        thin
    },
    liveedge/.style={
        draw=black,
        ->,
        line width=0.6pt
    }
]

% Root
\node[root] (r) at (6,0) {$r=8$};

% First level
\foreach \j/\x/\val/\sty in {
0/0/{1}/term,
1/1/{\ge 2}/cont,
2/2/{1}/term,
3/3/{1}/term,
4/4/{1}/term,
5/5/{1}/term,
6/6/{1}/term,
7/7/{1}/term,
8/8/{1}/term,
9/9/{1}/term,
10/10/{1}/term,
11/11/{1}/term,
12/12/{1}/term}
{
  \node[\sty] (v1\j) at (\x,-2) {$\val$};
  \node[anchor=south east,font=\tiny]
    at ([xshift=-0.4mm,yshift=0.4mm]v1\j.north west) {$v_{1,\j}$};
}

\foreach \x in {v10,v12,v13,v14,v15,v16,v17,v18,v19,v110,v111,v112}
    \draw[deadedge] (r.south) -- (\x.north);
\draw[liveedge] (r.south) -- (v11.north);

% Second level
\foreach \j/\x/\val/\sty in {
0/0/{2}/term,
1/1/{\ge 3}/cont,
2/2/{2}/term,
3/3/{2}/term,
4/4/{2}/term,
5/5/{2}/term,
6/6/{2}/term,
7/7/{2}/term,
8/8/{2}/term,
9/9/{2}/term,
10/10/{2}/term,
11/11/{2}/term,
12/12/{2}/term}
{
  \node[\sty] (v2\j) at (\x,-4.3) {$\val$};
  \node[anchor=south east,font=\tiny]
    at ([xshift=-0.4mm,yshift=0.4mm]v2\j.north west) {$v_{2,\j}$};
}

\foreach \x in {v20,v22,v23,v24,v25,v26,v27,v28,v29,v210,v211,v212}
    \draw[deadedge] (v11.south) -- (\x.north);
\draw[liveedge] (v11.south) -- (v21.north);

% Third level: no labels
\foreach \j/\x/\val/\sty in {
0/0/{3}/term,
1/1/{3}/term,
2/2/{3}/term,
3/3/{3}/term,
4/4/{3}/term,
5/5/{3}/term,
6/6/{3}/term,
7/7/{3}/term,
8/8/{\ge 4}/cont,
9/9/{3}/term,
10/10/{3}/term,
11/11/{3}/term,
12/12/{3}/term}
{
  \node[\sty] (v3\j) at (\x,-6.6) {$\val$};
}

\foreach \x in {v30,v31,v32,v33,v34,v35,v36,v37,v39,v310,v311,v312}
    \draw[deadedge] (v21.south) -- (\x.north);
\draw[liveedge] (v21.south) -- (v38.north);

% Fourth level: no labels
\foreach \j/\x/\val/\sty in {
0/0/{4}/term,
1/1/{4}/term,
2/2/{4}/term,
3/3/{4}/term,
4/4/{4}/term,
5/5/{4}/term,
6/6/{4}/term,
7/7/{4}/term,
8/8/{4}/term,
9/9/{4}/term,
10/10/{4}/term,
11/11/{4}/term,
12/12/{\ge 5}/cont}
{
  \node[\sty] (v4\j) at (\x,-8.9) {$\val$};
}

\foreach \x in {v40,v41,v42,v43,v44,v45,v46,v47,v48,v49,v410,v411}
    \draw[deadedge] (v38.south) -- (\x.north);
\draw[liveedge] (v38.south) -- (v412.north);

\node[anchor=north west,font=\small]
    at ([xshift=1mm,yshift=-1mm]current bounding box.north west)
    {\hspace{2cm}$T_{13}(8)$};
\end{tikzpicture}%
}
\caption{The beginning of the valuation tree rooted at $r=8$. Each box
records the common value of $\nu_{13}(e_n)$ on the corresponding family,
or a lower bound when the branch continues. The highlighted branch marks
the unique descendant that survives at each level.}
\label{fig:valuation-tree-8}
\end{figure}
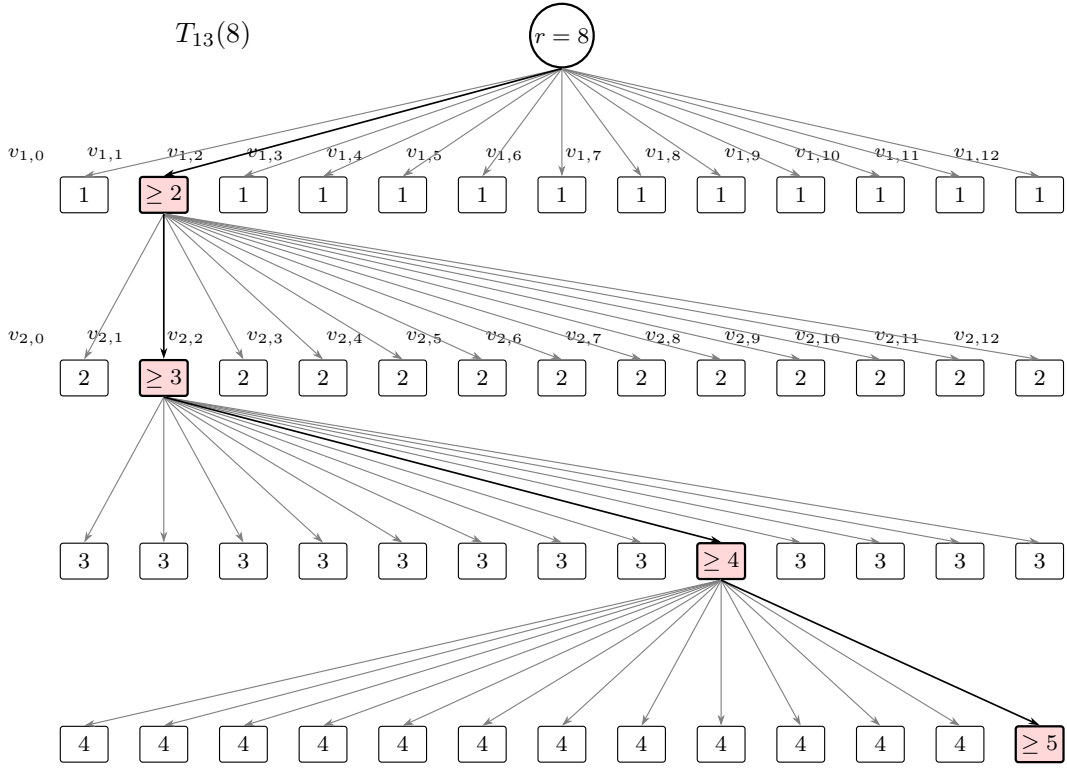

The refinement process may be continued inductively. At each step, a surviving congruence class modulo $2\cdot 13^k$ is refined into $13$ congruence classes modulo $2\cdot 13^{k+1}$. A vertex becomes terminal once the exact valuation is determined. In contrast, a vertex labeled by a lower bound, for example $\nu_{13}(e_n)\geq k+1$, must be refined at the next level. Figure~\ref{fig:valuation-tree-8} shows the first several levels of the tree $T_{13}(8)$.

The tree $T_{13}(8)$ is not an isolated phenomenon. The class
$I_{13}(12)$ behaves in exactly the same way. At the first refinement,
twelve of its children have valuation $1$, while a unique child has
valuation at least $2$:
\[
12\pmod{26}
\quad\longrightarrow\quad
38\pmod{338}.
\]
The classes $I_{13}(14)$ and $I_{13}(18)$ exhibit the same behavior,
although their surviving children occur at a different position:
\[
14\pmod{26}
\quad\longrightarrow\quad
300\pmod{338},
\]
and
\[
18\pmod{26}
\quad\longrightarrow\quad
304\pmod{338}.
\]
Finally, the remaining class $I_{13}(26)$ also has a unique surviving
child:
\[
26\pmod{26}
\quad\longrightarrow\quad
338\pmod{338}.
\]
There is also a striking symmetry among the four nonzero residue classes.
Indeed,
\[
8+18=26,
\qquad
12+14=26,
\]
and the corresponding surviving residues at the next level satisfy
\[
34+304=338,
\qquad
38+300=338.
\]
Thus the pairs $T_{13}(8)$ and $T_{13}(18)$, and $T_{13}(12)$ and
$T_{13}(14)$, already appear to be related by reflection.

The example $p=13$ therefore suggests a remarkably simple branching pattern. Each relevant residue class modulo $26$ splits into $13$ children at the next level. Of these children, twelve have their valuations completely determined and become terminal, while exactly one continues to a higher valuation. This same structure is visible at the successive levels of $T_{13}(8)$, as shown in Figure~\ref{fig:valuation-tree-8}.

The example $p=13$ leads us to ask whether this branching pattern persists for other primes and whether the unique surviving child, when it exists, can be determined without computing the often very large terms $e_n$ directly.

Motivated by these questions, we now formalize the refinement process described above through the notion of a $p$-adic valuation tree. For an odd prime $p$, a residue $r$, and $k\geq 0$, let
\[
I_{p,k}(r)
=
\{n\geq 1:n\equiv r\pmod{2p^{k+1}}\}.
\]
Thus, the classes considered in our example are $I_{13}(8)=I_{13,0}(8)$, 
and its first refinement consists of the $13$ classes
\[
I_{13,1}(8+26j),
\qquad 0\leq j\leq 12.
\]
More generally, every class $I_{p,k}(r)$ decomposes as
\[
I_{p,k}(r)
=
\bigcup_{j=0}^{p-1}
I_{p,k+1}\bigl(r+2p^{k+1}j\bigr).
\]
This elementary refinement is the underlying structure of the valuation
trees.

\begin{definition} 
Let $p$ be an odd prime and let $r$ satisfy $p\mid e_r$. The \emph{$p$-adic valuation tree rooted at $r$}, denoted by $T_p(r)$, is defined recursively as follows.

The root represents the class $I_{p,0}(r)$, and each vertex at level $k$ represents a congruence class $I_{p,k}(a)$. Such a vertex is called \emph{terminal} if $\nu_p(e_n)$ is constant as $n\in I_{p,k}(a)$; otherwise, it is called \emph{nonterminal}.

Each nonterminal vertex representing $I_{p,k}(a)$ has $p$ children, corresponding to the classes 
\[ 
I_{p,k+1}\bigl(a+2p^{k+1}j\bigr), \qquad 0\leq j\leq p-1. 
\] 
Terminal vertices are not refined further.
\end{definition}

The tree constructed in Figure~\ref{fig:valuation-tree-8} is precisely
$T_{13}(8)$. Its root represents $I_{13,0}(8)$, and its first generation
consists of the classes
\[
I_{13,1}(8+26j),
\qquad 0\leq j\leq12.
\]
The unique nonterminal vertex at this level corresponds to $j=1$, that
is, to the class $I_{13,1}(34)$.

Among the roots modulo $2p$, the class $2p\pmod{2p}$ plays a
distinguished role. Following Givens and Moll, a root
$r\in\{1,\ldots,2p-1\}$ for which $p\mid e_r$ is called an
\emph{internal root}, and the corresponding tree $T_p(r)$ is called an
\emph{internal branch}. Thus, for
$p=13$, the roots $8$, $12$, $14$, and $18$ are internal, whereas $26=2p$ is the distinguished root.

\begin{definition}\label{def:simple-tree}
An internal branch  $T_p(r)$ is called \emph{simple} if every level
contains exactly one nonterminal vertex.
\end{definition}

In particular, an internal root that is already terminal does not define
a simple branch. The computations above suggest that the internal branches
$T_{13}(8)$, $T_{13}(12)$, $T_{13}(14)$, and $T_{13}(18)$ are simple.
Givens and Moll introduced a difference function that not only proves
this fact, but also explains why exactly one child survives at each
level.

\subsection{The lifting criterion of Givens and Moll}

We now recall the lifting criterion of Givens and Moll, which provides a sufficient condition for an internal branch to be simple and determines the unique child on which the valuation may increase.

For $n\geq1$ and $j\geq0$, define
\begin{equation}\label{eq:Pnj}
P(n,j)=
\begin{cases}
1, & j=0,\\[1mm]
\displaystyle\prod_{\ell=1}^{j}(n^2-\ell^2), & j\geq1.
\end{cases}
\end{equation}
Givens and Moll showed that the even-indexed terms of the FiboLaplace sequence satisfy
\begin{equation*}
e_{2n}
=
n\sum_{j=0}^{n-1}P(n,j).
\end{equation*}
Accordingly, we define
\begin{equation}\label{eq:A-def}
A(n)
=
\frac{e_{2n}}{n}
=
\sum_{j=0}^{n-1}P(n,j).
\end{equation}
Let $p$ be an odd prime. For $r\geq1$, define
\begin{equation}\label{eq:Delta-GM}
\Delta_p(n,j;r)
\equiv
\frac{P(n+p^r,j)-P(n,j)}{p^r}
\pmod p
\end{equation}
and
\begin{equation}\label{eq:delta-GM}
\delta_p(n)
\equiv
\sum_{j=0}^{p}\Delta_p(n,j;r)
\pmod p.
\end{equation}
The quantity $\delta_p(n)$ is independent of $r$. The following result combines \cite[Lemma~4.8 and Theorem~4.12]{Givens}.

\begin{theorem}[Givens and Moll]\label{thm:GM-lifting}
Let $p$ be an odd prime.

\begin{enumerate}
\item For all integers $n,r\geq1$ and $j\geq0$,
\begin{equation}\label{eq:GM-even-lifting}
A(n+jp^r)-A(n)
\equiv
j\delta_p(n)p^r
\pmod{p^{r+1}}.
\end{equation}

\item Let $a\in\{1,\ldots,p-1\}$ satisfy $p\mid e_{2a}$. If
\[
\delta_p(a)\not\equiv0\pmod p,
\]
then the internal branch $T_p(2a)$ is simple; that is, every level
contains exactly one nonterminal vertex.
\end{enumerate}
\end{theorem}

To explain how the lifting congruence determines the children of a vertex, suppose that a nonterminal vertex at level $r-1$ is represented, in half-index notation, by a class $n\pmod{p^r}$, where $n\equiv a\pmod p$. Write
\[
A(n)\equiv cp^r\pmod{p^{r+1}}.
\]
Since $a\in\{1,\ldots,p-1\}$, we have $p\nmid n$. Therefore, $A(n)$ and $e_{2n}=nA(n)$ have the same $p$-adic valuation.
 
The $p$ children of this vertex are represented by
\[
n+jp^r,
\qquad
j=0,1,\ldots,p-1.
\]
Dividing \eqref{eq:GM-even-lifting} by $p^r$ gives
\begin{equation}\label{eq:linear-children}
\frac{A(n+jp^r)}{p^r}
\equiv
c+j\delta_p(n)
\pmod p.
\end{equation}
If $\delta_p(n)\not\equiv0\pmod p$, the linear congruence
\[
c+j\delta_p(n)\equiv0\pmod p
\]
has a unique solution $j$ modulo $p$. Consequently, exactly one child has valuation at least $r+1$. On each of the remaining $p-1$ child classes, the valuation is constantly equal to $r$, by periodicity modulo $p^{r+1}$, these children are terminal. Repeating this argument at each successive nonterminal vertex shows that $T_p(2a)$ is simple.

We now return to the example $p=13$. The internal root $8$ corresponds
to the half-index $a=4$. Since
\[
A(4)=\frac{e_8}{4}=1456=13\cdot112,
\]
we have $c\equiv112\equiv8\pmod{13}$. Moreover,
\[
\delta_{13}(4)\equiv5\pmod{13}.
\]
The first surviving child is therefore determined by
\[
8+5j\equiv0\pmod{13}.
\]
The unique solution is $j=1$, giving the full index $2(4+13)=34$. Thus the branch begins
\[
8\pmod{26}
\longrightarrow
34\pmod{338},
\]
as shown in Figure~\ref{fig:valuation-tree-8}.

The remaining three internal roots correspond to the half-indices
$6$, $7$, and $9$, for which
\[
\delta_{13}(6)\equiv11,\qquad
\delta_{13}(7)\equiv2,\qquad
\delta_{13}(9)\equiv8
\pmod{13}.
\]
Together with $\delta_{13}(4)\equiv5\pmod{13}$, these values are all
nonzero. Therefore, Theorem~\ref{thm:GM-lifting} proves that the four
internal branches rooted at  $8$, $12$,  $14$, and $18$ are simple.

Givens and Moll also obtained a parallel lifting criterion for the odd
subsequence; see \cite[Example~4.14]{Givens}. Set
\begin{equation}\label{eq:B-def}
B(n)=e_{2n+1}
=
1+\sum_{j=1}^{n}n(n+j)P(n,j-1),
\qquad n\geq0.
\end{equation}
For $r\geq1$, define
\begin{equation}\label{eq:delta-prime}
\delta'_p(n)
\equiv
\frac{B(n+p^r)-B(n)}{p^r}
\pmod p.
\end{equation}
The quotient is an integer, and its residue modulo $p$ is independent
of $r$.

\begin{theorem}[Givens and Moll]\label{thm:GM-odd}
Let $p$ be an odd prime.

\begin{enumerate}
\item For all integers $n\geq0$, $r\geq1$, and $j\geq0$,
\begin{equation}\label{eq:GM-odd-lifting}
B(n+jp^r)-B(n)
\equiv
j\delta'_p(n)p^r
\pmod{p^{r+1}}.
\end{equation}

\item Let $a\in\{0,1,\ldots,p-1\}$ satisfy $p\mid e_{2a+1}$. If
\[
\delta'_p(a)\not\equiv0\pmod p,
\]
then the internal branch $T_p(2a+1)$ is simple.
\end{enumerate}
\end{theorem}

The preceding results led Givens and Moll to ask whether every even
internal branch must be simple.

\begin{OpenP}[Givens and Moll]\label{prob:GM}
Let $p$ be an odd prime and let $a\in\{1,\ldots,p-1\}$ satisfy
$p\mid e_{2a}$. Must the internal branch $T_p(2a)$ be simple?
\end{OpenP}

Theorem~\ref{thm:GM-lifting} gives the sufficient condition
\[
\delta_p(a)\not\equiv0\pmod p
\quad\Longrightarrow\quad
T_p(2a)\text{ is simple}.
\]
This leads naturally to the stronger question of whether this
nonvanishing condition holds automatically at every even internal root;
that is, whether
\[
p\mid e_{2a}
\quad\Longrightarrow\quad
\delta_p(a)\not\equiv0\pmod p.
\]

It is important to distinguish these two questions. Givens and Moll
proved that the nonvanishing of $\delta_p(a)$ implies simplicity, but
they did not prove the converse. Thus, finding an internal root for which $\delta_p(a)\equiv0\pmod p$ would disprove the stronger nonvanishing statement, but would not by
itself show that the corresponding branch is not simple.

Our counterexample shows that both assertions fail simultaneously. For
a certain even internal root, we have $\delta_p(a)\equiv0\pmod p$, and the root itself is terminal. Consequently, the corresponding
internal branch is not simple.

\section{A Terminal Even Internal Branch} \label{sec:counterexample}

The vanishing of $\delta_p(a)$ does not, by itself, imply that the corresponding internal branch is not simple. It only shows that the sufficient criterion of Givens and Moll does not apply. When the initial valuation is exactly one, however, the lifting congruence gives a stronger conclusion: the root is terminal. We use this observation to show that $T_{971}(388)$ is not simple, thereby giving a negative answer to the question of Givens and Moll in the even case. The analogous odd counterexample is more immediate: the valuation formula of Givens and Moll for $p=3$ shows directly that the branch $T_3(3)$ is terminal.

\begin{proposition}\label{prop:terminal-criterion}
Let $p$ be an odd prime and let $a\in\{1,\ldots,p-1\}$. If $\nu_p(e_{2a})=1$ and $\delta_p(a)\equiv0\pmod p$, then the root of the internal branch $T_p(2a)$ is terminal. More precisely,
\begin{equation}\label{eq:terminal-observation}
\nu_p\bigl(e_{2(a+jp)}\bigr)=1
\qquad\text{for every }j\geq0.
\end{equation}
\end{proposition}
\begin{proof}
Since $p\nmid a$, we have $\nu_p(A(a))=\nu_p(e_{2a})=1$. Taking $n=a$ and $r=1$ in \eqref{eq:GM-even-lifting} gives
\[
A(a+jp)-A(a)
\equiv
j\delta_p(a)p
\equiv0
\pmod{p^2}.
\]
Hence
\[
A(a+jp)\equiv A(a)\pmod{p^2}
\]
for every $j\geq0$. Since $A(a)$ is divisible by $p$ but not by $p^2$,
the same is true of $A(a+jp)$. Moreover, $p\nmid a+jp$, and therefore $\nu_p\bigl(e_{2(a+jp)}\bigr)=1$ for every $j\geq0$. Thus the valuation is constant on the root class $2a\pmod{2p}$, so the
root is terminal.
\end{proof}

\subsection{A terminal branch for \texorpdfstring{$p=971$}{p=971}} 
We now apply the preceding criterion to construct an even internal branch whose root is terminal. This provides the counterexample needed to answer the question of Givens and Moll negatively.

The following counterexample was discovered through a computer search implemented in Mathematica. 

The prime $p=971$ provides precisely such a situation. For $a=194$, exact modular computations give 
\begin{equation}\label{eq:971-data} 
e_{388} \equiv 473848 = 971\cdot488 \pmod{971^2} 
\end{equation} 
and 
\begin{equation}\label{eq:971-delta} 
\delta_{971}(194)\equiv0\pmod{971}. 
\end{equation} 
Both congruences can be verified using modular arithmetic. The first is obtained by iterating the recurrence \eqref{enRecurrenceRelation} modulo $971^2$, while the second follows from the recurrence for $\Delta_p(n,j;r)$ established by Givens and Moll. Since $971\nmid488$, equation \eqref{eq:971-data} implies that 
\[ 
\nu_{971}(e_{388})=1. 
\] 
Therefore, Proposition~\ref{prop:terminal-criterion}, together with \eqref{eq:971-delta}, yields \[ \nu_{971}(e_{388+1942j})=1 \qquad\text{for every }j\geq0. \] Thus, the root of the internal branch $T_{971}(388)$ is terminal, so $T_{971}(388)$ is not simple. This provides a negative answer to the question of Givens and Moll.

The natural odd analogue is already false for $p=3$. Indeed, Givens
and Moll proved in \cite[Theorem~3.6]{Givens} that
$\nu_3(e_m)=1$ whenever $m\equiv3\pmod6$. Hence the odd internal root $3$ is terminal, and therefore
$T_3(3)$ is not simple.

The vanishing of the corresponding lifting parameter can also be seen
directly. Since
\[
B(1)=e_3=3
\qquad\text{and}\qquad
B(4)=e_9\equiv3\pmod9,
\]
we have
\[
\delta'_3(1)
\equiv
\frac{B(4)-B(1)}{3}
\equiv0
\pmod3.
\]
Thus both the tree question and the stronger nonvanishing question have
negative answers in the odd case.

\section{A polynomial reformulation of the Givens--Moll criterion}
\label{sec:polynomial-reformulation}

The counterexamples in the preceding section show that the nonvanishing conditions in the Givens--Moll lifting criteria are not automatic at internal roots. We now interpret these conditions algebraically by associating polynomials $\mathcal P_p(x)$ and $\mathcal Q_p(x)$ with the even and odd subsequences. We show that internal roots correspond to roots of these polynomials modulo $p$, while the lifting parameters correspond to their derivatives. This yields simple-root criteria for the simplicity of even and odd internal branches. We also analyze multiple roots and establish, when the initial $p$-adic valuation is one, a partial converse in the even case.

\subsection{The even polynomial reformulation}
We begin with the even subsequence. We introduce a polynomial $\mathcal P_p(x)$ whose nonzero zeros modulo $p$ correspond precisely to the even internal root vertices of the valuation tree. We then show that the derivative of $\mathcal P_p(x)$ determines the even lifting parameter $\delta_p(a)$. This interpretation yields a criterion for the simplicity of even internal branches in terms of simple zeros of $\mathcal P_p(x)$ and clarifies the behavior when the corresponding zero is multiple.

For each $k\geq0$, the expression $P(n,k)$ in \eqref{eq:Pnj} is a polynomial in $n$. We then define
 
\begin{equation}\label{eq:root-polynomial}
\mathcal P_p(x)
=
\sum_{k=0}^{p-1}P(x,k)
=
1+\sum_{k=1}^{p-1}\prod_{\ell=1}^{k}(x^2-\ell^2).
\end{equation}

For $a\in\{1,\ldots,p-1\}$, the terms  $P(a,k)$ vanish for $k\geq a$. Hence $A(a)=\mathcal P_p(a)$, and therefore $e_{2a}=a\mathcal P_p(a)$. Since $p\nmid a$, it follows that
\[
p\mid e_{2a}
\quad\Longleftrightarrow\quad
\mathcal P_p(a)\equiv0\pmod p.
\]
Thus the even internal roots $2a$ correspond precisely to the nonzero
roots $a$ of $\mathcal P_p(x)$ modulo $p$.

We shall use the following elementary congruence.
\begin{lemma}\label{lem:first-order-polynomial-congruence}
Let $f(x)\in\mathbb Z[x]$, let $p$ be a prime, and let $r\geq1$.
For all integers $a$ and $j$,
\begin{equation}\label{eq:first-order-polynomial-congruence}
f(a+jp^r)
\equiv
f(a)+jp^rf'(a)
\pmod{p^{r+1}}.
\end{equation}
\end{lemma}
\begin{proof}
Write $f(x)=\sum_{m=0}^d c_mx^m$. By the binomial theorem,
\[
(a+h)^m =\sum_{i=0}^m\binom{m}{i}a^{m-i}h^{i}=
a^m+ma^{m-1}h+h^2Q_m(a,h),
\]
for some polynomial $Q_m(a,h)$ with integer coefficients. Hence
\[
f(a+h)
=
f(a)+hf'(a)+h^2Q(a,h),
\]
where $Q(a,h)$ is again a polynomial in $a$ and $h$ with integer
coefficients. Taking $h=jp^r$, the last term is
divisible by $p^{2r}$, and since $r\geq1$, we have $2r\geq r+1$.
Therefore
\[
f(a+jp^r)
\equiv
f(a)+jp^rf'(a)
\pmod{p^{r+1}}. \qedhere
\]
\end{proof}

The next result gives the desired interpretation of the quantity
$\delta_p(a)$.

\begin{theorem}\label{thm:delta-derivative}
Let $p$ be an odd prime. For every $a\in\{1,\ldots,p-1\}$,
\begin{equation}\label{eq:delta-derivative}
\delta_p(a)
\equiv
\mathcal P_p'(a)
\pmod p.
\end{equation}
\end{theorem}

\begin{proof}
Applying Lemma~\ref{lem:first-order-polynomial-congruence} to
$f(x)=P(x,k)$, with $j=1$, gives
\[
P(a+p^r,k)-P(a,k)
\equiv
p^r\frac{\partial}{\partial x}P(x,k)\bigg|_{x=a}
\pmod{p^{r+1}}.
\]
Dividing by $p^r$ and using the definition of
$\Delta_p(a,k;r)$ in \eqref{eq:Delta-GM}, we obtain
\[
\Delta_p(a,k;r)
\equiv
\frac{\partial}{\partial x}P(x,k)\bigg|_{x=a}
\pmod p.
\]
Therefore, by \eqref{eq:delta-GM},
\[
\delta_p(a)
\equiv
\sum_{k=0}^{p}
\frac{\partial}{\partial x}P(x,k)\bigg|_{x=a}
\pmod p
\equiv
\mathcal P_p'(a)
+
\frac{\partial}{\partial x}P(x,p)\bigg|_{x=a}
\pmod p.
\]
It remains to show that the last term vanishes modulo $p$. By the product rule,
\[
\frac{\partial}{\partial x}P(x,p)
=
\sum_{j=1}^{p}
2x\prod_{\substack{1\leq \ell\leq p\\ \ell\neq j}}
(x^2-\ell^2).
\]
Evaluating at $x=a$, every term with $j\neq a$ vanishes, since its
product contains the factor $a^2-a^2$. Hence
\[
\frac{\partial}{\partial x}P(x,p)\bigg|_{x=a}
=
2a\prod_{\substack{1\leq\ell\leq p\\ \ell\neq a}}
(a^2-\ell^2).
\]
Since $p$ is odd, $p-a\neq a$, and therefore the factor corresponding
to $\ell=p-a$ occurs in this product. Moreover,
$a^2-(p-a)^2=p(2a-p)$, so the product is divisible by $p$. Thus
\[
\frac{\partial}{\partial x}P(x,p)\bigg|_{x=a}
\equiv0\pmod p,
\]
and \eqref{eq:delta-derivative} follows.
\end{proof}

It follows that, whenever $2a$ is an even internal root,  $\delta_p(a)\not\equiv0\pmod p$ if and only if $a$ is a simple root of $\mathcal P_p(x)$ modulo $p$.
Thus the stronger nonvanishing question arising from the work of Givens
and Moll has the following exact algebraic reformulation:
\[
\begin{split}
&p\mid e_{2a}
\Longrightarrow
\delta_p(a)\not\equiv0\pmod p
\quad\text{for every }a\in\{1,\ldots,p-1\}
\\
&\hspace{25mm}\Longleftrightarrow\quad
\text{every nonzero root of $\mathcal P_p(x)$ modulo $p$ is simple}.
\end{split}
\]

This gives the following corrected positive statement.

\begin{theorem}\label{thm:simple-root-criterion}
Let $p$ be an odd prime, and let $a\in\{1,\ldots,p-1\}$ satisfy
$p\mid e_{2a}$. If $a$ is a simple root of $\mathcal P_p(x)$ modulo $p$,
then the internal branch $T_p(2a)$ is simple.
More precisely, set $c\equiv \frac{A(a)}{p}\pmod p$. Then the unique nonterminal child of the root corresponds to the
half-index $a+jp$, where $j\in\{0,\ldots,p-1\}$ is the unique solution of
\begin{equation}\label{eq:first-surviving-child-polynomial}
c+j\mathcal P_p'(a)\equiv0\pmod p.
\end{equation}
\end{theorem}

\begin{proof}
Since $a$ is a simple root, $\mathcal P_p'(a)\not\equiv0\pmod p$. By Theorem~\ref{thm:delta-derivative}, $\delta_p(a)\not\equiv0\pmod p$, and Theorem~\ref{thm:GM-lifting} therefore implies that
$T_p(2a)$ is simple. 

To determine the first surviving child, take $n=a$ and $r=1$ in
\eqref{eq:GM-even-lifting}. Dividing by $p$ gives
\[
\frac{A(a+jp)}{p}
\equiv
c+j\delta_p(a)
\equiv
c+j\mathcal P_p'(a)
\pmod p.
\]
Since $\mathcal P_p'(a)\not\equiv0\pmod p$, the congruence
\eqref{eq:first-surviving-child-polynomial} has a unique solution
$j\in\{0,\ldots,p-1\}$.
\end{proof}

\begin{corollary}\label{cor:all-even-branches-simple}
If every nonzero root of $\mathcal P_p(x)$ modulo $p$ is simple, then
every even internal branch of the valuation tree is simple.
\end{corollary}

%\begin{proof}
%Let $T_p(2a)$ be an even internal branch. Then
%$a\in\{1,\ldots,p-1\}$ and $p\mid e_{2a}$, so $\mathcal P_p(a)\equiv0\pmod p$. Thus $a$ is a nonzero root of $\mathcal P_p(x)$ modulo $p$, and by
%hypothesis it is simple. The result now follows from
%Theorem~\ref{thm:simple-root-criterion}.
%\end{proof}

\begin{remark}\label{rem:Hensel}
The unique-solution argument in
\eqref{eq:first-surviving-child-polynomial} is the simple-root case of
Hensel lifting. No form of Hensel's lemma is needed here, however; the
lifting congruence \eqref{eq:GM-even-lifting} gives the required
conclusion directly.
\end{remark}

It remains to understand what the same congruence says when $a$ is a multiple root of $\mathcal P_p(x)$ modulo $p$. Suppose that
\[
A(a)\equiv cp\pmod{p^2}
\qquad\text{and}\qquad
\mathcal P_p'(a)\equiv0\pmod p.
\]
Then
\[
\frac{A(a+jp)}{p}\equiv c\pmod p
\]
for every $j\in\{0,\ldots,p-1\}$. If $c\not\equiv0\pmod p$, then all $p$ children have valuation exactly $1$. Hence the root is terminal. If $c\equiv0\pmod p$, then all $p$
children have valuation at least $2$. In this case, the first-order
congruence does not distinguish among them, and higher-order information
is needed to determine which children have valuation exactly $2$ and
which, if any, have higher valuation. Thus multiple roots are precisely the cases in which the first-order
Givens--Moll criterion no longer determines a unique surviving child.

\subsection{The odd polynomial reformulation} We now develop the corresponding polynomial reformulation for the odd subsequence. We introduce a polynomial $\mathcal Q_p(x)$ whose nonzero roots modulo $p$ correspond precisely to the odd internal root vertices of the valuation tree. We then show that the derivative of $\mathcal Q_p(x)$ determines the odd lifting parameter $\delta'_p(a)$. As in the even case, this interpretation yields a simple-root criterion for the simplicity of odd internal branches and clarifies the behavior when the corresponding polynomial root is multiple.

Set $Q(x,0)=1$, and for $k\geq1$ define
\[
Q(x,k)
=
\prod_{\ell=0}^{k-1}(x-\ell)(x+\ell+1).
\]
Thus \eqref{eq:B-def} can be written as
\[
B(n)=\sum_{k=0}^{n}Q(n,k).
\]
Define
\begin{equation}\label{eq:Qcal-def}
\mathcal Q_p(x)
=
\sum_{k=0}^{p-1}Q(x,k).
\end{equation}

For $a\in\{0,\ldots,p-1\}$, the terms $Q(a,k)$ vanish for
$k\geq a+1$. Hence $B(a)=\mathcal Q_p(a)$, and therefore
\begin{equation}\label{eq:odd-root-characterization}
p\mid e_{2a+1}
\quad\Longleftrightarrow\quad
\mathcal Q_p(a)\equiv0\pmod p.
\end{equation}
Thus the odd internal roots $2a+1$ correspond precisely to the roots
$a$ of $\mathcal Q_p(x)$ modulo $p$.
 
The odd lifting parameter also has a derivative interpretation.
\begin{proposition}\label{prop:odd-derivative}
For every $a\in\{0,\ldots,p-1\}$,
\[
\delta'_p(a)
\equiv
\mathcal Q_p'(a)
\pmod p.
\]
\end{proposition}

The proof is entirely analogous to the even case: one applies
Lemma~\ref{lem:first-order-polynomial-congruence} to each polynomial
$Q(x,k)$ and then uses the definition of $\delta'_p(a)$.

\begin{theorem}\label{thm:odd-simple-root}
Let $p$ be an odd prime, and let $a\in\{0,\ldots,p-1\}$ satisfy
$p\mid e_{2a+1}$. If $a$ is a simple root of $\mathcal Q_p(x)$
modulo $p$, then the odd internal branch $T_p(2a+1)$ is simple.
\end{theorem}

\begin{proof}
Since $a$ is a simple root of $\mathcal Q_p(x)$ modulo $p$, we have $\mathcal Q_p'(a)\not\equiv0\pmod p$. By Proposition~\ref{prop:odd-derivative},
$\delta'_p(a)\equiv \mathcal Q_p'(a)\not\equiv0\pmod p$, and the conclusion follows from Theorem~\ref{thm:GM-odd}.
\end{proof}

\begin{corollary}\label{cor:all-odd-branches-simple}
If every root of $\mathcal Q_p(x)$ modulo $p$ is simple, then every odd
internal branch of the valuation tree is simple.
\end{corollary}

The behavior at a multiple root is also parallel to the even case.
Suppose that
\[
B(a)\equiv cp\pmod{p^2}
\qquad\text{and}\qquad
\mathcal Q_p'(a)\equiv0\pmod p.
\]
By Proposition~\ref{prop:odd-derivative} and
\eqref{eq:GM-odd-lifting},
\[
\frac{B(a+jp)}{p}\equiv c\pmod p
\]
for every $j\geq0$. If $c\not\equiv0\pmod p$, then
\[
\nu_p(e_{2(a+jp)+1})=1
\qquad\text{for every }j\geq0,
\]
so the odd internal root $2a+1$ is terminal. If
$c\equiv0\pmod p$, all $p$ children acquire a second factor of $p$,
and higher-order information is needed to distinguish them.

Finally, the odd polynomial satisfies the symmetry
$\mathcal Q_p(-x-1)=\mathcal Q_p(x)$. Differentiating
 and evaluating at the fixed point $a=(p-1)/2$ modulo $p$ gives
\[
\mathcal Q_p'\left(\frac{p-1}{2}\right)\equiv0\pmod p.
\]
Therefore, whenever the central odd index $p$ is an internal root, the
corresponding root of $\mathcal Q_p(x)$ is necessarily multiple. This
is precisely the phenomenon that occurs for $p=3$.

\subsection{A partial converse}
The simple-root criterion gives a sufficient condition for an even internal branch to be simple. Under the additional assumption that the initial $p$-adic valuation is exactly one, the converse also holds. Thus, in this case, simplicity of the branch is equivalent to both the nonvanishing of the Givens--Moll lifting parameter and the simplicity of the corresponding polynomial root of $\mathcal P_p(x)$ modulo $p$.

\begin{theorem}
\label{thm:converse-simple-root}
Let $p$ be an odd prime, and let $a\in\{1,\ldots,p-1\}$ satisfy $\nu_p(e_{2a})=1$. Then the following conditions are equivalent:
\begin{enumerate}
\item the even internal branch $T_p(2a)$ is simple;
\item $\delta_p(a)\not\equiv0\pmod p$;
\item $a$ is a simple root of $\mathcal P_p(x)$ modulo $p$.
\end{enumerate}
\end{theorem}

\begin{proof}
The equivalence between (2) and (3) follows from
Theorem~\ref{thm:delta-derivative}, since
\[
\delta_p(a)\equiv\mathcal P_p'(a)\pmod p,
\]
and $a$ is already a root of $\mathcal P_p(x)$ modulo $p$.

The implication $(3)\Rightarrow(1)$ is precisely
Theorem~\ref{thm:simple-root-criterion}.

It remains to prove $(1)\Rightarrow(3)$. Suppose that $T_p(2a)$ is
simple and that $\mathcal P_p'(a)\equiv0\pmod p$. Since $p\nmid a$ and $e_{2a}=a\mathcal P_p(a)$, the assumption
$\nu_p(e_{2a})=1$ implies
\[
\frac{\mathcal P_p(a)}{p}\not\equiv0\pmod p.
\]
By the multiple-root analysis above, this forces the root of
$T_p(2a)$ to be terminal, contradicting the assumption that the branch
is simple. Hence $\mathcal P_p'(a)\not\equiv0\pmod p$, so $a$ is a simple root of $\mathcal P_p(x)$ modulo $p$.
\end{proof}

\section{A recursive $p$-adic lifting procedure}
\label{sec:padic-lifting}

The results of the preceding section provide criteria for determining when an internal branch is simple. We now show that, for a simple branch, the lifting congruence can also be used recursively to follow its unique nonterminal vertex. We begin with the case $p=13$, which illustrates the procedure and shows how the exact $13$-adic valuation of a term can be determined from a sequence of successively refined residue classes.

\subsection{The case $p=13$}
\label{subsec:p13-lifting}

We return to the case $p=13$ considered earlier. Modulo $13$, the polynomial $\mathcal P_{13}(x)$ factors as 

\[ 
\mathcal P_{13}(x) \equiv{} (x^2-3)(x^2+3)(x^6+x^4+5x^2-3) (x^{14}+6x^{10}+2x^6-5x^4-3x^2+1) \pmod{13}. 
\] The last two factors have no roots in $\mathbb F_{13}$. Moreover, the roots of $x^2-3$ are $4$ and $9$, while the roots of $x^2+3$ are $6$ and $7$. Hence the nonzero roots of $\mathcal P_{13}(x)$ modulo $13$ are 
\[ 
4,\quad 6,\quad 7,\quad 9. 
\] 
They correspond, respectively, to the even internal root vertices $8$, $12$, $14$, and $18$, and hence to the internal branches 
\[ 
T_{13}(8),\quad T_{13}(12),\quad T_{13}(14),\quad T_{13}(18). 
\]

To determine the first surviving child of each branch, we use
\eqref{eq:first-surviving-child-polynomial}. Since
$A(a)=\mathcal P_{13}(a)$ for $1\leq a\leq12$, set
\[
c_a
\equiv
\frac{\mathcal P_{13}(a)}{13}
\pmod{13},
\qquad
d_a
\equiv
\mathcal P_{13}'(a)
\pmod{13}.
\]

By \eqref{eq:delta-derivative}, $d_a\equiv\delta_{13}(a)\pmod{13}$.
Thus the first surviving child is determined by the unique solution
$j_a\in\{0,\ldots,12\}$ of
\[
c_a+j_ad_a\equiv0\pmod{13},
\]
or equivalently,
\[
j_a\equiv-c_ad_a^{-1}\pmod{13}.
\]
Equation~\eqref{eq:root-polynomial} yields the following data:
\[
\begin{array}{c|cccc}
a&4&6&7&9\\
\hline
\displaystyle
c_a=\frac{\mathcal P_{13}(a)}{13}\pmod{13}
&8&2&4&3\\[2mm]
\displaystyle
d_a=\mathcal P_{13}'(a)\pmod{13}
&5&11&2&8\\
\displaystyle
j_a\equiv-c_ad_a^{-1}\pmod{13}
&1&1&11&11
\end{array}
\]
In particular, $d_a\not\equiv0\pmod{13}$ for all four roots. Hence
they are simple, and Theorem~\ref{thm:simple-root-criterion} shows that
the corresponding four internal branches are simple.

For instance, for $a=4$, equation
\eqref{eq:first-surviving-child-polynomial} becomes
\[
8+5j\equiv0\pmod{13}.
\]
Its unique solution is $j=1$. Hence the surviving half-index modulo
$13^2$ is $4+13=17$. Since the corresponding full index is twice the half-index, the first
step of the branch is
\[
8\pmod{26}
\longrightarrow
34\pmod{338}.
\]
Applying the same calculation to the other three roots gives
\[
\begin{aligned}
12\pmod{26}
&\longrightarrow38\pmod{338},\\
14\pmod{26}
&\longrightarrow300\pmod{338},\\
18\pmod{26}
&\longrightarrow304\pmod{338}.
\end{aligned}
\]
In each case, the remaining twelve children are terminal with valuation
exactly $1$.

The procedure can now be repeated at the unique surviving child. We
illustrate this again for the branch rooted at $8$. Its first surviving
half-index is $17$. Since
\[
\frac{A(17)}{13^2}\equiv8\pmod{13},
\]
and the slope along this branch is $\mathcal P_{13}'(4)\equiv5\pmod{13}$, the next surviving digit is determined by
\[
8+5j\equiv0\pmod{13}.
\]
Thus $j=1$, and the next surviving half-index is $17+13^2=186$. At the following step,
\[
\frac{A(186)}{13^3}\equiv12\pmod{13},
\]
so the corresponding congruence is
\[
12+5j\equiv0\pmod{13}.
\]
Its unique solution is $j=8$, giving the next surviving half-index $186+8\cdot13^3=17762$.

Therefore, in terms of the full index, the branch begins
\[
8\pmod{26}
\longrightarrow
34\pmod{338}
\longrightarrow
372\pmod{4394}
\longrightarrow
35524\pmod{57122}
\longrightarrow\cdots.
\]
At each level, the surviving residue modulo $2\cdot13^r$ splits into thirteen possible lifts modulo $2\cdot13^{r+1}$, and the lifting
congruence selects exactly one of them.

There are no odd internal branches for $p=13$. Indeed, a direct
computation shows that $\mathcal Q_{13}(x)$ has no roots in
$\mathbb F_{13}$. Hence, by \eqref{eq:odd-root-characterization},
no odd residue class gives rise to an internal branch.

It remains to consider the distinguished even class. We have
\[
\mathcal P_{13}(0)\equiv1\pmod{13}.
\]
If $13\mid m$, then \eqref{eq:GM-even-lifting}, with $r=1$, gives
\[
A(m)\equiv A(13)\pmod{13}.
\]
Moreover, from the definitions of $A$ and $\mathcal P_{13}$,
\[
A(13)
=
\sum_{k=0}^{12}P(13,k)
\equiv
\sum_{k=0}^{12}P(0,k)
=
\mathcal P_{13}(0)
\pmod{13}.
\]
Hence $A(m)\equiv1\pmod{13}$. Thus $A(m)$ is a $13$-adic unit, and since $e_{2m}=mA(m)$, 
we obtain
\[
\nu_{13}(e_{2m})=\nu_{13}(m).
\]
Consequently, at level $r$, the distinguished even class has the
surviving residue $0\pmod{2\cdot13^r}$.

\begin{figure}[ht!]
\centering
\begin{tikzpicture}[ scale=0.72,
    transform shape,
    x=1cm,
    y=1cm,
    every node/.style={font=\small},
    root/.style={circle,fill=black,inner sep=2.6pt},
    terminal/.style={circle,fill=gray!55,draw=gray!70,inner sep=1.3pt},
    bbranch/.style={circle,fill=blue!75,inner sep=2.2pt},
    gbranch/.style={circle,fill=green!60!black,inner sep=2.2pt},
    obranch/.style={circle,fill=orange!90,inner sep=2.2pt},
    pbranch/.style={circle,fill=violet!80,inner sep=2.2pt},
    terminaledge/.style={gray!65,line width=.35pt},
    bedge/.style={blue!70,line width=.9pt},
    gedge/.style={green!50!black,line width=.9pt},
    oedge/.style={orange!85!red,line width=.9pt},
    pedge/.style={violet!80,line width=.9pt}
]

\usetikzlibrary{calc,positioning}

%------------------------------------------------
% Root
%------------------------------------------------
\node[root] (R) at (0,0) {};
\node[above=2pt of R] {$\mathrm{root}$};

%------------------------------------------------
% Level 1: classes modulo 26
% 22 terminal + 4 nonterminal
%------------------------------------------------
\node[bbranch] (B1) at (-7.2,-1.8) {};
\node[gbranch] (G1) at (-2.4,-1.8) {};
\node[obranch] (O1) at ( 2.4,-1.8) {};
\node[pbranch] (P1) at ( 7.2,-1.8) {};

\draw[bedge] (R)--(B1);
\draw[gedge] (R)--(G1);
\draw[oedge] (R)--(O1);
\draw[pedge] (R)--(P1);

\node[blue!75,left=3pt of B1] {$8$};
\node[green!50!black,left=3pt of G1] {$12$};
\node[orange!85!red,right=3pt of O1] {$14$};
\node[violet!80,right=3pt of P1] {$18$};

% 22 terminal children at level 1
\foreach \x in {
-10.2,-9.4,-8.6,-7.8,-6.4,-5.6,-4.8,-4.0,-3.2,-1.6,-0.8,
 0.8, 1.6, 3.2, 4.0, 4.8, 5.6, 6.4, 7.8, 8.6, 9.4,10.2}
{
    \draw[terminaledge] (R)--(\x,-1.8);
    \node[terminal] at (\x,-1.8) {};
}

%------------------------------------------------
% Surviving vertices
%------------------------------------------------
\node[bbranch] (B2) at (-6.4,-3.7) {};
\node[bbranch] (B3) at (-7.0,-5.6) {};
\node[bbranch] (B4) at (-6.3,-7.5) {};

\node[gbranch] (G2) at (-1.7,-3.7) {};
\node[gbranch] (G3) at (-2.3,-5.6) {};
\node[gbranch] (G4) at (-1.6,-7.5) {};

\node[obranch] (O2) at ( 1.7,-3.7) {};
\node[obranch] (O3) at ( 2.3,-5.6) {};
\node[obranch] (O4) at ( 1.6,-7.5) {};

\node[pbranch] (P2) at ( 6.4,-3.7) {};
\node[pbranch] (P3) at ( 7.0,-5.6) {};
\node[pbranch] (P4) at ( 6.3,-7.5) {};

%------------------------------------------------
% Macro: 12 terminal children around a parent
%------------------------------------------------
\newcommand{\terminalfan}[2]{%
    \foreach \dx in {-2.0,-1.6,-1.2,-0.8,-0.4,-0.15,0.15,0.4,0.8,1.2,1.6,2.0}
    {
       \node[terminal] at ($(#1)+(\dx,#2)$) {};
       \draw[terminaledge] (#1)--($(#1)+(\dx,#2)$);
    }
}

%------------------------------------------------
% Blue branch: 8 -> 34 -> 372 -> 21598
%------------------------------------------------
\terminalfan{B1}{-1.45}
\draw[bedge] (B1)--(B2);
\node[blue!75,right=4pt of B2] {$34$};

\terminalfan{B2}{-1.45}
\draw[bedge] (B2)--(B3);
\node[blue!75,right=4pt of B3] {$372$};

\terminalfan{B3}{-1.45}
\draw[bedge] (B3)--(B4);
\node[blue!75,right=4pt of B4] {$21598$};

\terminalfan{B4}{-1.45}

%------------------------------------------------
% Green branch: 12 -> 38 -> 638 -> 25726
%------------------------------------------------
\terminalfan{G1}{-1.45}
\draw[gedge] (G1)--(G2);
\node[green!50!black,right=4pt of G2] {$38$};

\terminalfan{G2}{-1.45}
\draw[gedge] (G2)--(G3);
\node[green!50!black,right=4pt of G3] {$638$};

\terminalfan{G3}{-1.45}
\draw[gedge] (G3)--(G4);
\node[green!50!black,right=4pt of G4] {$25726$};

\terminalfan{G4}{-1.45}

%------------------------------------------------
% Orange branch: 14 -> 300 -> 3756 -> 31396
%------------------------------------------------
\terminalfan{O1}{-1.45}
\draw[oedge] (O1)--(O2);
\node[orange!85!red,left=4pt of O2] {$300$};

\terminalfan{O2}{-1.45}
\draw[oedge] (O2)--(O3);
\node[orange!85!red,left=4pt of O3] {$3756$};

\terminalfan{O3}{-1.45}
\draw[oedge] (O3)--(O4);
\node[orange!85!red,left=4pt of O4] {$31396$};

\terminalfan{O4}{-1.45}

%------------------------------------------------
% Violet branch: 18 -> 304 -> 4022 -> 35524
%------------------------------------------------
\terminalfan{P1}{-1.45}
\draw[pedge] (P1)--(P2);
\node[violet!80,left=4pt of P2] {$304$};

\terminalfan{P2}{-1.45}
\draw[pedge] (P2)--(P3);
\node[violet!80,left=4pt of P3] {$4022$};

\terminalfan{P3}{-1.45}
\draw[pedge] (P3)--(P4);
\node[violet!80,left=4pt of P4] {$35524$};

\terminalfan{P4}{-1.45}

\end{tikzpicture}
\caption{The four simple internal branches for $p=13$.}
\label{fig:p13-fractal-tree}
\end{figure}
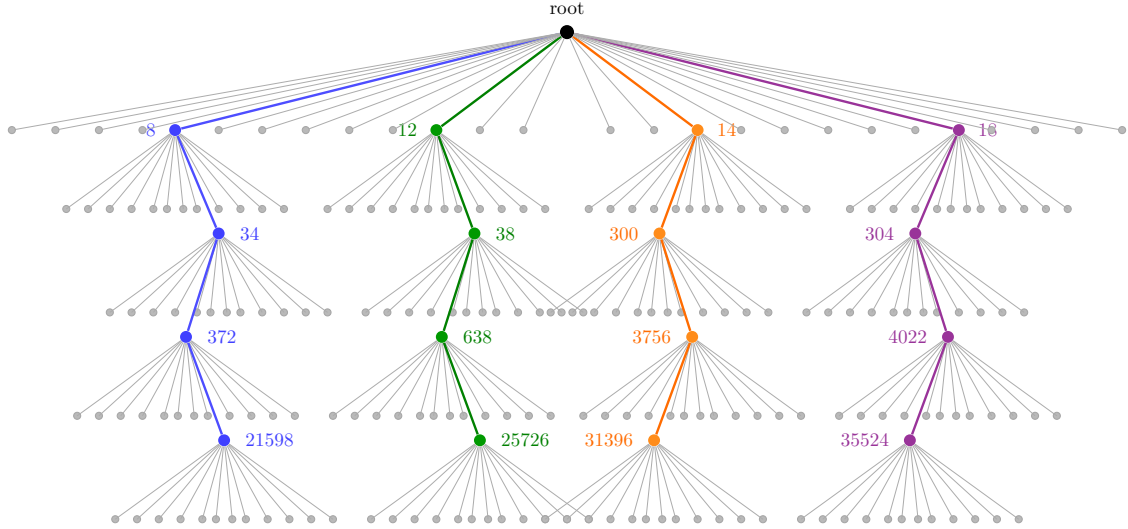

Figure~\ref{fig:p13-fractal-tree} summarizes the lifting behavior of the
four internal branches. The colored paths record the unique surviving
full-index residue at each successive level, while the gray vertices
represent terminal classes. For clarity, the distinguished even branch
corresponding to the residue $0$ is not displayed; as observed above,
this branch survives independently at every level.

Let $\mathcal R_r$ denote the set of full-index residue classes modulo
$2\cdot13^r$ for which $13^r\mid e_n$. The first four sets are
\[
\begin{array}{c|c|c}
r&2\cdot13^r&\mathcal R_r\\
\hline
1&26&
\{0,8,12,14,18\}\\
2&338&
\{0,34,38,300,304\}\\
3&4394&
\{0,372,638,3756,4022\}\\
4&57122&
\{0,21598,25726,31396,35524\}.
\end{array}
\]
Thus, for every $r\geq1$,
\[
13^r\mid e_n
\quad\Longleftrightarrow\quad
n\bmod(2\cdot13^r)\in\mathcal R_r.
\]

For example,
\[
372\equiv8\pmod{26},
\qquad
372\equiv34\pmod{338},
\qquad
372\equiv372\pmod{4394},
\]
whereas $372\notin\mathcal R_4$. Therefore
\[
\nu_{13}(e_{372})=3.
\]

This example suggests a general procedure. Rather than constructing the
entire valuation tree, one may follow only the unique surviving residue
of each simple branch. We now formalize this lifting process.

\subsection{The general lifting procedure}
\label{subsec:general-lifting-procedure}

The preceding example suggests how to proceed for an arbitrary odd
prime $p$. Recall that
\[
A(m)=\frac{e_{2m}}{m},
\qquad
B(m)=e_{2m+1},
\]
and define
\[
E_p=
\left\{
a\in\{1,\ldots,p-1\}:
\mathcal P_p(a)\equiv0\pmod p
\right\},
\]
and
\[
O_p=
\left\{
b\in\{0,\ldots,p-1\}:
\mathcal Q_p(b)\equiv0\pmod p
\right\}.
\]
By the results of the preceding section, the elements of $E_p$
correspond to the nonzero even internal branches $T_p(2a)$, while the
elements of $O_p$ correspond to the odd internal branches
$T_p(2b+1)$. The distinguished even class $0\pmod{2p}$ will be treated
separately.

We now describe how to follow the successive surviving residues of any
simple branch.

\begin{theorem}
\label{thm:recursive-simple-lifting}
Let $p$ be an odd prime. If $a\in E_p$ is a simple root of $\mathcal P_p(x)$ modulo $p$, then
there exists a unique sequence
\[
u_1^{(a)},u_2^{(a)},u_3^{(a)},\ldots,
\qquad
0\leq u_r^{(a)}<p^r,
\]
such that
\[
u_1^{(a)}=a,
\qquad
u_{r+1}^{(a)}\equiv u_r^{(a)}\pmod{p^r},
\qquad
p^r\mid A\bigl(u_r^{(a)}\bigr)
\]
for every $r\geq1$.

Similarly, if $b\in O_p$ is a simple root of $\mathcal Q_p(x)$ modulo
$p$, then there exists a unique sequence
\[
v_1^{(b)},v_2^{(b)},v_3^{(b)},\ldots,
\qquad
0\leq v_r^{(b)}<p^r,
\]
such that
\[
v_1^{(b)}=b,
\qquad
v_{r+1}^{(b)}\equiv v_r^{(b)}\pmod{p^r},
\qquad
p^r\mid B\bigl(v_r^{(b)}\bigr)
\]
for every $r\geq1$.
\end{theorem}

\begin{proof}
We prove the even case; the odd case is analogous. Suppose that
$u_r^{(a)}$ has been determined and set
\[
c_r^{(a)}
\equiv
\frac{A\bigl(u_r^{(a)}\bigr)}{p^r}
\pmod p.
\]
The $p$ possible lifts are $u_r^{(a)}+jp^r$, for $j=0,\ldots,p-1$. By \eqref{eq:GM-even-lifting},
\[
\frac{A\bigl(u_r^{(a)}+jp^r\bigr)}{p^r}
\equiv
c_r^{(a)}
+
j\delta_p\bigl(u_r^{(a)}\bigr)
\pmod p.
\]
Since $u_r^{(a)}\equiv a\pmod p$, Theorem~\ref{thm:delta-derivative}
gives
\[
\delta_p\bigl(u_r^{(a)}\bigr)
\equiv
\mathcal P_p'\bigl(u_r^{(a)}\bigr)
\equiv
\mathcal P_p'(a)
\pmod p.
\]
Therefore
\[
\frac{A\bigl(u_r^{(a)}+jp^r\bigr)}{p^r}
\equiv
c_r^{(a)}+j\mathcal P_p'(a)
\pmod p.
\]
Since $\mathcal P_p'(a)\not\equiv0\pmod p$, there is a unique
$j\in\{0,\ldots,p-1\}$ for which the right-hand side vanishes.
This determines the unique lift divisible by $p^{r+1}$ and completes
the recursion.
\end{proof}

The proof is constructive: once the surviving residue at a given level is known, the lifting congruence uniquely determines the next base-$p$ digit and, hence, the surviving residue at the next level. The odd case is entirely analogous. Consequently, along every simple branch, the exact $p$-adic valuation of a term can be determined one base-$p$ digit at a time.

\section*{Acknowledgments}

R. Fl\'orez, was supported by the Citadel Foundation. R. Higuita was supported by Universidad Pedag\'ogica y Tecnol\'ogica de Colombia under Grant number SGI 4141. 
J. L. Ramírez was partially supported by Universidad Nacional de Colombia, Project No. 66041. 

\section*{Statement on computational tools}
The counterexample in Section~\ref{sec:counterexample} was
discovered through a search implemented in Mathematica. ChatGPT
was used to assist with the Mathematica code and with the
presentation of the manuscript. All computational results were
verified by the authors.

\medskip

\noindent MSC2020: 11B39, 33B10

\end{document}